\documentclass[a4paper,12pt]{article}
\usepackage{amsmath,amsfonts,amssymb}

\newtheorem{theorem}{Theorem}
\newtheorem{proposition}[theorem]{Proposition}

\def\proof{\noindent{\it Proof.}\ }

\def\be{\begin{equation}}
\def\dsp{\displaystyle }

\def\a{\alpha}
\def\b{\beta}
\def\c{\gamma}
\def\d{\delta}
\def\l{\lambda}

\def\ss{\sigma}

\def\A{{\mathbb A}}
\def\B{{\mathbb B}}
\def\N{{\mathbb N}}
\def\C{{\mathbb C}}
\def\Z{{\mathbb Z}}

\title{ About division by 1}
\author{ Alain Lascoux\footnote{text written during the Conference
{\it Applications of the Macdonald Polynomials},
at the Newton Institute in April 2001. }} 

\begin{document}
\date{}
\maketitle

\begin{abstract}
The Euclidean division of two formal series in one variable 
produces a sequence of series that we obtain explicitly, 
remarking that the case
where one of the two initial series is 1 is sufficiently generic. 
As an application, we define a Wronskian of symmetric functions.
\end{abstract}

The Euclidean division of two polynomials $P(z)$, $Q(z)$, in one variable
$z$, of consecutive degrees, produces a sequence of linear factors 
(the successive quotients), and a sequence of successive remainders,
both families being symmetric functions in the roots of $P$ and $Q$ 
separately.

Euclidean division can also be applied to formal series in $z$, but it never
stops in the generic case, leaving time enough to observe the law 
of the coefficients appearing in the process.

Moreover, since the quotient of two formal series is also a formal series,
it does not make much difference if we
 suppose that one of the two initial series
is 1. This renders the division of series simpler than that of polynomials;
in fact the latter could be obtained from the former. 

By {\it formal series} we mean a unitary series 
$$ f(z) = 1 +c_1 z +c_2 z^2 + \cdots  \ .$$
We shall moreover formally factorize it
$$ f(z) =\sigma_z(\A) := \prod\nolimits_{a\in \A} (1-za)^{-1} = 
  \sum_{i=0}^\infty z^i\, S_i(\A) \ ,$$
the {\it alphabet} $\A$ being supposed to be an infinite set of indeterminates,
or of complex numbers, the coefficients $S_i(\A)$ being called
the {\it complete functions} in $\A$.

Given two series, dividing $f_{-1}(z)=\sigma_z(\A)$ by 
$f_0(z)= \sigma_z(\B)$ means finding the unique coefficients
$\a,\b$ such that
\begin{equation} 
\bigl(\ss_z(\A) -(1+\a z)\, \ss_z(\B)  \bigr)\, \frac{1}{\b}\, z^{-2} 
\end{equation}
is a unitary series $f_1(z)= \ss_z(\C)$.

Dividing in turn $f_0(z)$ by $f_1(z)$, one obtains $f_2(z)$,
and iterating one gets from the intial pair $(f_{-1},\, f_0)$ 
an infinite sequence of series $f_{-1},\, f_0,\, f_1,\, f_2,\, f_3, \ldots $

However, all the equations
\be f_{k-1}(z)= (1+\a_k\, z)\, f_k(z) + \b_k\,z^2\, f_{k+1}(z) \label{eq:2} 
\end{equation}
can be divided by $f_0(z) = \ss_z(\B)$.
If the $k$-th remainder for the pair $(\ss_z(\A),\, \ss_z(\B))$ is
$\ss_z(\C)$, then the $k$-th remainder for the pair 
$\bigl(\ss_z(\A)/\ss_z(\B)\, ,\, 1  \bigr)$ is 
$$ \ss_z(\C -\B) = \sum_{i=0}^\infty z^i\, S_i(\C-\B) 
 := \frac{\prod_{b\in \B} (1-zb)}{\prod_{c\in \C} (1-zc) }\ .  $$
Indeed, already at the first step, one sees that the equation
\be \ss_z(\A) = \bigl(1+zS_1(\A-\B) \bigr)\, \ss_z(\B) + 
   S_2(\A-\B)\, z^2\, f_1(z)  \label{eq:3}
\end{equation}
is equivalent to
\be \ss_z(\A-\B) = \Bigl(1+zS_1(\A-\B) \Bigr)\, 1 +
   S_2(\A-\B)\, z^2\, f_1(z)/\ss_z(\B)   \label{eq:4}
\end{equation}
since the two coefficients are functions of $\A-\B$. 

Of course equation (\ref{eq:4}) expands into 
\be \ss_z(\A-\B) = \bigl(1+zS_1(\A-\B) \bigr) +
   S_2(\A-\B)\, z^2 \sum_i z^i  \frac{S_{2+i}(\A-\B)}{S_2(\A-\B)} \label{eq:5}
\end{equation}

We shall see that the other remainders can be as easily written. We just need
to recall the definition of a Schur function $S_\l(\A-\B)$, where 
$\l =[\l_1\geq \l_2\geq \cdots \geq \l_\ell\geq 0]$ is a partition~:
$$ S_\l(\A-\B) = \det\Bigl( S_{\l_i+j-i}(\A-\B) \Bigr)_{1\leq i,j \leq \ell}
\ . $$
Then
\begin{theorem}
 The $k$-th remainder in the Euclidean division of 
$\ss_z(\A)$ by $1$ is
\be
 f_k(z) = S_{(k+1)^k}(\A)^{-1}\, \sum_{i=0}^\infty
 z^i\, S_{k+1+i,(k+1)^{k-1}}(\A) \ .
\end{equation}
The $k$-th remainder of the division of $\ss_z(\A)$
by $\ss_z(\B)$ is
\be
 f_k(z) = S_{(k+1)^k}(\A-\B)^{-1}\, \ss_z(\B)\,  \sum_{i=0}^\infty
 z^i\, S_{k+1+i,(k+1)^{k-1}}(\A-\B) \ .  \label{eq:7}
\end{equation}

\end{theorem}

\proof Merging all equations (\ref{eq:2}) together, one can characterize
in the first case 
$f_k(z)$ as the unique series such that there exists scalars
$\a_1^k,\ldots, \a_{k-1}^k$, $\b_1^k,\ldots, \b_k^k$, $\gamma^k$ :
\be
 z^{2k}\gamma^k\, f_k(z) = 
(1+\a_1^k z +\cdots +\a_{k-1}^k z^{k-1})\, \ss_z(\A)
 - (1+ \b_1^k z +\cdots + \b_k^k z^k)\, 1 \ . \label{eq:8}
\end{equation}
Defining
$$ \ss_y(\A \pm \frac{1}{z}) = \ss_y(\A) \ss_y(\pm \frac{1}{z}))
= \ss_y(\A) (1-y/z) ^{\mp 1} \ ,$$
I claim that 
\be 
z^{k-1} S_{(k+1)^{k-1}}(\A -\frac{1}{z})\, \ss_z(\A) 
+ (-z)^k S_{k^k}(\A +\frac{1}{z})  \label{eq:9}
\end{equation}
is equal to $z^{2k} S_{(k+1)^k}(\A)$ modulo terms of higher
degree in $z$.

Indeed, 
$z^{k-1} S_{(k+1)^{k-1}}(\A -1/z)\, \ss_z(\A)$ can be written as
the $k\times k$ determinant
$$ \left|
\begin{matrix}
   z^{k-1}\ss_z(\A) & \cdots & \ss_z(\A) \\
  S_{k}(\A) & \cdots & S_{2k-1}(\A)  \\
  \vdots      &        & \vdots      \\
 S_2(\A) & \cdots & S_{k+1}(\A) \\
\end{matrix}\right| $$
because, subtracting from each column, except the first,
 $1/z$ times the preceding column 
and using $S_j(\A-1/z)= S_j(\A) -S_{j-1}(\A)/z $, $j\in \Z$, 
the determinant factorizes
  into $z^{k-1}\ss_z(\A)\, S_{(k+1)^{k-1}}(\A -1/z)$.

The coefficients of $z^{k+1},\ldots, z^{2k-1}$ are the functions
$S_{2,(k+1)^{k-1}}(\A), \ldots,$ 
$S_{k,(k+1)^{k-1}}(\A)$
which are zero, having two identical rows in their determinantal
expression.

The coefficients of $z^0,\ldots, z^k$ are
$$ S_{1-k,(k+1)^{k-1}}(\A) , \ldots,\, S_{1,(k+1)^{k-1}}(\A)$$
which, after permuting their first row with the others,
 are recognized to be
$$ S_{k^{k-1},0}(\A), \ldots, S_{k^{k-1},k}(\A) \ ,$$
up to the sign $(-1)^{k-1}$.

These terms sum up to $(-z)^k S_{k^k}(\A +1/z) $, 
and the required series $f_k(z)$ is obtained by dividing
by $(-1)^{k-1} S_{k^{k-1}}(\A)$.  \hfill QED

\medskip
The first defining equations for the $f_i$'s are,
writing $S_j$ for $S_j(\A)$ :
$$ z^2S_2\, f_1 = \ss_z(\A) -(1+zS_1) = \ss_z(\A) -zS_1(\A+\frac{1}{z}) \ ,$$
\begin{multline*}
 -z^4 S_{33}/S_2\, f_2= (1-zS_3/S_2)\ss_z(A) 
  - (1+z S_{21}/S_2 +z^2S_{22}/S_2) = \\
 = -S_3(\A-\frac{1}{z})\ss_z(\A)/S_2 -S_{22}(\A+\frac{1}{z}) /S_2 \ ,
\end{multline*}
\begin{multline*}
z^6S_{444}/S_{33}\, f_3 = 
 (1-zS_{43}/S_{33}+z^2 S_{44}/S_{33})\, \ss_z(\A) 
   -(1+z S_{331}/S_{33} +z^2S_{332}/S_{33} \\ +z^3 S_{333}/S_{33})  
  = z^2 S_{44}(\A-\frac{1}{z}) \, \ss_z(\A)/S_{33} -z^3 
S_{333}(\A+\frac{1}{z}) /S_{33} 
\end{multline*}

Equation (\ref{eq:9}) can be understood as giving the Pad\'e approximant
of degree $[k,k-1]$ of the series $\ss_z(\A)$~:
\be
\ss_z(\A) = (-1)^{k-1} \frac{z\, S_{k^k}(\A+1/z)}{S_{(k+1)^{k-1}}(\A-1/z)}
  + z^{2k} \frac{S_{(k+1)^k}(\A)}{z^{k-1}S_{(k+1)^{k-1}}(\A-1/z)}\, f_k(z)
\end{equation} 

\medskip Sylvester treated the Euclidean division by a different method,
using summations on subsets of roots (cf. [LP]). 
We shall for our part interpret now
the Euclidean division as producing a sequence
of alphabets from a given one.
This time, it is more convenient to divide 1 by a series,
and put $f_{-1}=1=\ss_z(0)$, $f_0=\ss_z(\A)=\ss_z(\A^0), 
\, \ldots ,\, f_k =\ss_z(\A^k),
\ldots $.

From (\ref{eq:7}) one has 
$$ \ss_z(\A^k -\A^0)= \
S_{(k+1)^k}(0-\A^0)^{-1}\,  \sum_{i=0}^\infty
 z^i\, S_{k+1+i,(k+1)^{k-1}}(0-\A^0)  $$
that is,
\be
S_{k^{k+1}}(\A)\, \ss_z(\A^k) = \ss_z(\A)
 \sum_{i=0}^\infty z^i\, S_{k+1+i,(k+1)^{k-1}}(-\A^0) \label{eq:12} 
\end{equation}

\begin{proposition}   The successive remainders $\ss_z(\A^k)$ 
in the division of 1 by
$\ss_z(\A)$ satisfy 
\be 
 S_{k^{k+1}}(\A)\, \ss_z(\A^k) =
\sum_{i=0}^\infty z^i\, S_{k+i, k^k}(\A)   \label{eq:13}
\end{equation}
\end{proposition}

\proof Instead of having recourse to determinants, and using relations between
minors, let us invoke symmetrizing operators.  Suppose the cardinality 
of $\A$ to be finite,
$\A=\{a_1,\ldots, a_N \}$, before letting it tend to infinity.
Let $\pi_\omega$ be the {\it maximal isobaric divided difference},
that is, the operator such that the image of 
$a^\l:= a_1^{\l_1}\cdots a_N^{\l_N}$ is $S_\l(\A)$,
for any $\l\in \N^N$.  The  series on the right of (\ref{eq:13})
is the image under $\pi_\omega$ of 
$$ (a_1\cdots a_{k+1})^k /(1-za_1) $$
and also of
 $$ (a_1\cdots a_{k+1})^k /(1-za_1)\cdots (1-za_{k+1}) \ , $$
which can be written 
$$ \frac{(a_1\cdots a_{k+1})^k\, (1-za_{k+2})\cdots (1-za_{N})}
{(1-za_1) \cdots (1-za_N)} \ .$$
The denominator is symmetrical in $a_1,\ldots, a_N$,
and thus commutes with $\pi_\omega$. 
As for the numerator, the only monomials giving a non-zero contribution
have exponents 
$$ \underbrace{k,\ldots, k}_{k+1} \underbrace{1,\ldots, 1}_{i} 
 \underbrace{0,\ldots, 0}_{N-k-1-i}
\ , \ 0\leq i\leq N-k-1 \ . $$
Therefore
$$  \sum_{i=0}^\infty S_{k+i, k^k}(\A) =
\sum (-z)^i S_{k^{k+1},1^i}(\A) \, \ss_z(\A)  \ ,$$
which gives (\ref{eq:12})  since
$S_\l(\A) = (-1)^{|\l|} S_{\l^\sim} (-\A)$, where 
$\l^\sim$ is the conjugate to the partition $\l$.   
\hfill QED

The sequences $\A, \A^1, \A^2, \ldots$ have been considered by B.Leclerc,
to whom  the following notion  of a Wronskian of 
complete symmetric functions is due 
(there are more general Wronskians associated to any set of
symmetric functions, and any alphabet).

Let $n$ be an integer, $k_1,\ldots, k_n\in \N$, and $\A$ be an alphabet.
Then the {\it Wronskian} $W( S_{k_1},\ldots,\, S_{k_n}\, ;\, \A)$
is the determinant
$$ \det \Bigl|S_{k_j-i}(\A_i)  \Bigr|_{0\leq i\leq n-1, 1\leq j\leq n} \ , $$
where $\A^0:=\A,\, \A^1,\, \A^2 \ldots$ is the sequence of alphabets
obtained in the Euclidean division of 1 by $\ss_z(\A)$.

As an application of his study of relations between minors [Bl],
B. Leclerc  obtained in an unpublished note~: 

\begin{theorem}
Let $n$ be an integer, $K=[k_1,\ldots, k_n]\in \N^n$, $\A$ be an alphabet.
Then 
\be
 W( S_{k_1},\ldots,\, S_{k_n}\, ;\, \A) =
S_{ K+ [n-1,\ldots,0]}(\A) / S_{(n-1)^n}(\A)  \ . \label{eq:14}
\end{equation}
\end{theorem}

\proof To simplify the notation, we shall take $n=4$, 
$K[\a,\b,\c,\d]$. Then the Wronskian, in terms of functions of
$\A$ only, is, according to (\ref{eq:13}), 
\be
\left|
 \begin{matrix}
  S_{\a000} & S_{\b000} & S_{\c 000} & S_{\d 000} \\
  S_{\a100} & S_{\b100} & S_{\c 100} & S_{\d 100} \\
  S_{\a220} & S_{\b220} & S_{\c 220} & S_{\d 220} \\
  S_{\a333} & S_{\b333} & S_{\c 333} & S_{\d 333} \\
\end{matrix} \right|
\end{equation} 
up to division of the second, third, last row by 
$S_{11}, S_{222}, S_{3333}$ respectively.

All the entries of this determinant are $4\times 4$ minors of
the matrix~:
\be 
\begin{matrix} {\scriptstyle column\ index} \\ \\ \\ \\ \\  
\end{matrix}
\left[
\begin{matrix}
    0  & 1   & 2   &3   &4   &5   &{\a}     &{\b}     &{\c}     &{\d} \\
   \noalign{\hrule}
   S_0 & S_1 & S_2 &S_3 &S_4 &S_5 &S_{\a+3} &S_{\b+3} &S_{\c+3} &S_{\d+3} \\
   .   & S_0 & S_1 & S_2 &S_3 &S_4 &S_{\a+2} &S_{\b+2} &S_{\c+2} &S_{\d+2} \\
   .   & . & S_0 & S_1 & S_2 &S_3   &S_{\a+1} &S_{\b+1} &S_{\c+1} &S_{\d+1} \\
   .   & . & .   & S_0 & S_1 & S_2  &S_{\a}  &S_{\b}  &S_{\c}   &S_{\d}  \\
\end{matrix}  \right]
\end{equation}

Designating a minor by the sequence  indexing its columns, the determinant
 to study is
$$ \left| 
\begin{matrix}
[012\a] & [012\b] &[012\c] &[012\d] \\\relax
[013\a] & [013\b] &[013\c] &[013\d] \\\relax
[034\a] & [034\b] &[034\c] &[034\d] \\\relax
[345\a] & [345\b] &[345\c] &[345\d] 
\end{matrix}
\right|
$$
and according to Bazin (cf. [Bl]), factorizes into
$$ [0123]\, [0134]\, [0345]\, [\a\b\c\d] \,  $$
that is, more explicitly, into the product
$$ S_{0000}(\A)\ S_{1100}(\A)\ S_{2220}(\A)\ S_{3333}(\A) \ . $$
Reintroducing the missing factor $\ S_{11}(\A), S_{222}(\A)$,
we obtain the theorem.  \hfill QED 

\smallskip
As we have said, Pad\'e approximants, formal orthogonal polynomials,
continued fraction expansion of formal series are all 
related to Euclidean division.

As a last example, we illustrate how to write a continued
fraction expansion of a series, which one can find in the work of 
 Wronski, Chebyschef or Stieljes.

\begin{proposition}  Given an alphabet $\A$, then

\smallskip\noindent
${1\over z}\, \ss_{1/z}(\A) = $ 
 $$
\qquad{}{1\over \dsp z+S_1(0-\A^0)+
{\strut S_2(0-\A^0) \over \dsp z+S_1(\A^0-\A^1) +
{\strut  S_2(\A^0-\A^1) \over \dsp  z+S_1(\A^1-\A^2) +
{\strut S_2(\A^1-\A^2)\over \dsp \ddots}}}}
$$
\end{proposition}

\proof The validity of such expansion amounts to the recursions
(\ref{eq:3}) 
$$ \ss_{1/z}(\A^k) = \Bigl( 1+ {1\over z} S_1(\A^{k-1} -\A^k)
\Bigr) \ss_{1/z}(\A^{k}) + S_2(\A^{k-1} -\A^k)\,
{1\over z^2} \ss_{1/z}(\A^{k+1})
$$
with which we began this text.

\medskip
\centerline{\bf References}
\medskip
\parindent 0mm

[Bl]\ B. Leclerc, 
{\it On identities satisfied by minor of a matrix}, Advances in Maths. 
 {\bf 100} (1993)101--132.

[LP]\  A. Lascoux, P. Pragacz,
  {Sylvester's formulas for Euclidean division and Schur multifunctions},
 Max Planck Institut Preprints (2001).

\vskip 15mm
{\obeylines
\hfill          C.N.R.S., Institut Gaspard Monge   \hfill
\hfill          Universit\'e de Marne-la-Vall\'ee,  \hfill
\hfill          5 Bd Descartes, Champs sur Marne,   \hfill
\hfill          77454 Marne La Vall\'ee Cedex 2 FRANCE  \hfill
\hfill \tt        Alain.Lascoux@univ-mlv.fr \hfill
\hfill          http://phalanstere.univ-mlv.fr/$\sim$al \hfill
}

\end{document}